\documentclass[english]{llncs}
\usepackage{amsmath}
\usepackage{amssymb}
\usepackage{bbm}
\usepackage{mathbbol}
\usepackage{mathrsfs}
\usepackage[all]{xy}
\newcounter{contatore}
\title{Lambda theories of effective lambda models}

\author{Chantal Berline\inst{1}
        \and
        Giulio Manzonetto\inst{1,2}
        \and
        Antonino Salibra\inst{2}\\
        chantal.berline@pps.jussieu.fr, \{gmanzone, salibra\}@dsi.unive.it
}

\institute{Laboratoire PPS, CNRS-Universit\'{e} Paris 7,\\
2, place Jussieu (case 7014), 75251 Paris Cedex 05, France
\and
Universit\`a Ca'Foscari di Venezia,\\
Dipartimento di Informatica
Via Torino 155,
30172 Venezia, Italy}

\begin{document}
% Macro di Giulio

\newcommand{\Tins}[1]{#1^{\cT-ins}}
\newcommand{\cardinality}[1]{card(#1)}% cardinality of a set 
\newcommand{\codepair}[2]{\hspace{+2pt}<\hspace{-2pt}#1,#2\hspace{-3pt}>}% code of a pair
\newcommand{\codepairmix}[2]{\hspace{+2pt}\ll\hspace{-2pt} #1,#2 \hspace{-2.5pt}\gg}% code of a pair N^* X N
\newcommand{\Visser}[1]{\mathscr{V}_{#1}}
\newcommand{\restr}{\upharpoonright}% restriction
\newcommand{\Kdn}[1]{\widehat{#1}}% restriction
\newcommand{\emptyrho}{\rho_\bot}% empty environment
\newcommand{\RE}{{\mathcal R}{\mathcal E}}% r.e. sets
\newcommand{\CORE}{{\mathcal R}{\mathcal E}^{co}}% co-r.e.
\newcommand{\Kleeneq}{\simeq}% Kleene equivalence
\newcommand{\Zero}{\textrm{Zero}}
\newcommand{\Pred}{\textrm{Pred}}
\newcommand{\Succ}{\textrm{Succ}}
\newcommand{\ecomp}{\zeta}
\newcommand{\st}{:}
\newcommand{\pow}[1]{\mathscr{P}(#1)}
\newcommand{\Intunder}[1]{\Lambda^o_{#1}}
\newcommand{\Prime}[1]{\cP(#1)}
\newcommand{\step}[2]{\varepsilon_{#1,#2}}
\newcommand{\D}{{\mathcal D}}
\newcommand{\E}{{\mathcal E}}
\newcommand{\std}[1]{\underline{#1}}
\newcommand{\cod}[1]{\ulcorner #1 \urcorner}
\newcommand{\vphi}{\varphi}
\newcommand{\partarrow}{\rightharpoonup}
\newcommand{\W}{{\cal W}}
\newcommand{\App}{\cF}
\newcommand{\Abs}{\lambda}
\newcommand{\compl}[1]{#1^c}
\newcommand{\domcomp}[1]{#1^{comp}}
\newcommand{\domdec}[1]{#1^{dec}}
\newcommand{\Intcomp}[1]{\left|  #1\right|^{comp}}
\newcommand{\Int}[1]{\left|  #1\right|}
\newcommand{\clup}[1]{#1\hspace{-5pt}\uparrow}
\newcommand{\cldn}[1]{#1\hspace{-5pt}\downarrow}
\newcommand{\inv}[1]{#1^{-}}
\newcommand{\partinv}[1]{#1^{-1}}
\newcommand{\img}[1]{#1^{+}}
\renewcommand{\sup}[1]{\bigsqcup #1}
\newcommand{\Unsolvable}{\mathscr{U}}
\newcommand{\cl}[1]{\overline{#1}}
\newcommand{\seq}[1]{\overline{#1}}
\newcommand{\const}[1]{\underline{#1}}
\newcommand{\intrho}[1]{|#1|^i_\rho}
\renewcommand{\bold}[1]{{\bf #1}}
\newcommand{\compel}[1]{{\bf {\cal K}}(#1)}
\newcommand{\Appr}[1]{\c App(\textrm{$#1$})}
\newcommand{\imp}{\Rightarrow}
\newcommand{\sqle}{\sqsubseteq}
\newcommand{\BTle}{\sqsubseteq_\cB}
\newcommand{\BTeq}{=}
\newcommand{\BTsup}{\bigsqcup}
\newcommand{\nat}{\mathbb{N}}
\newcommand{\BT}{\mathcal{BT}}
\newcommand{\Base}{\mathcal{O}}

% Macro di Nino
%\newtheorem{conjecture}{Conjecture}

%\newtheorem{theorem}{Theorem}[section]
%\newtheorem{notation}[theorem]{Notation}
%\newtheorem{corollary}[theorem]{Corollary}
%\newtheorem{proposition}[theorem]{Proposition}
%\newtheorem{lemma}[theorem]{Lemma}
%\newtheorem{definition}[theorem]{Definition}
%\newtheorem{remark}[theorem]{Remark}
%\newtheorem{claim}[theorem]{Claim}
%\newtheorem{subclaim}[theorem]{Subclaim}
%\newtheorem{example}[theorem]{Example}

%
\newcommand{\ot}{\overline{t}}
\newcommand{\os}{\overline{s}}
\newcommand{\sskV}{{\bf k}\sp {\bf V}}
\newcommand{\sssV}{{\bf s}\sp{\bf V}}
\newcommand{\glV}{{\gl^{\bf V}}}
\newcommand{\cdotV}{\cdot^{\bf V}}
\newcommand{\cterm}{combinatory term}
\newcommand{\bZd}{\mbox{\bf Zd\,}}
\newcommand{\Cr}{\mbox{\rm Cr\,}}
\newcommand{\ConA}{{\bf C}{\rm on}\hspace{.01in}{\bf A}}
\newcommand{\ConB}{{\bf C}{\rm on}\hspace{.01in}{\bf B}}
\newcommand{\ConAB}{{\bf C}{\rm on}\hspace{.01in}({\bf A}$\times${\bf B})}
\newcommand{\lvar}{$\gl$-variable}
\newcommand{\bI}{{\bf I }}
\newcommand{\TolA}{{\bf T}{\rm ol}\hspace{.01in}{\bf A}}
\newcommand{\FRA}{FR\hspace{.01in}{\bf A}}
\newcommand{\bnew}{\marginpar{\mbox{$\top$}}}
\newcommand{\enew}{\marginpar{\mbox{$\bot$}}}
\newcommand{\snew}{\marginpar{\mbox{$\leftarrow$}}}
\newcommand{\fnew}{\marginpar{\mbox{$\downarrow$}}}
\newcommand{\glxA}{\gl x^{\bA}}
\newcommand{\Zd}{\mbox{\rm Zd\,}}
% theorem-like commands
%

%\newtheorem{claim}{Claim}[section]
%\newtheorem{theorem}{Theorem}[section]
%\newtheorem{corollary}{Corollary}[section]
%\newtheorem{definition}{Definition}
%\newtheorem{definition}{Definition}[section]
%\newtheorem{proposition}{Proposition}[section]
%\newtheorem{lemma}{Lemma}[section]
%\newtheorem{example}{Example}[section]

\newtheorem{claimsub}{Claim}[subsection]
\newtheorem{theoremsub}{Theorem}[subsection]
\newtheorem{corollarysub}{Corollary}[subsection]
\newtheorem{definitionsub}{Definition}[subsection]
\newtheorem{propositionsub}{Proposition}[subsection]
\newtheorem{lemmasub}{Lemma}[subsection]
\newtheorem{examplesub}{Example}[subsection]

\newtheorem{claimsubsub}{Claim}[subsubsection]
\newtheorem{theoremsubsub}{Theorem}[subsubsection]
\newtheorem{corollarysubsub}{Corollary}[subsubsection]
\newtheorem{definitionsubsub}{Definition}[subsubsection]
\newtheorem{propositionsubsub}{Proposition}[subsubsection]
\newtheorem{lemmasubsub}{Lemma}[subsubsection]
\newtheorem{examplesubsub}{Example}[subsubsection]

\newcommand{\Over}[1]{\raisebox{-1.6mm}{{\scriptsize #1}}}
\newcommand{\SubOver}[1]{\raisebox{-0.6mm}{{\tiny #1}}}

\newcommand{\FOL}{\mbox{$L_{\omega,\omega}$}}

\newcommand{\Neg}{\mbox{$\neg$}}
\newcommand{\noteq}{\mbox{$\neq$}}
\newcommand{\fiproof}{\mbox{$(\Rightarrow)$}}
\newcommand{\ifproof}{\mbox{$(\Leftarrow)$}}
\newcommand{\of}{\mbox{$\in$}}
\newcommand{\notof}{\mbox{$\not\in$}}
\newcommand{\union}{\mbox{$\cup$}}
\newcommand{\intersection}{\mbox{$\cap$}}
\newcommand{\except}{\mbox{$-$}}
\newcommand{\cardof}[1]{\mbox{$|$}#1\mbox{$|$}}
\newcommand{\forevery}{\mbox{$\forall$}}
\newcommand{\therexists}{\mbox{$\exists$}}
\newcommand{\within}{\mbox{$\subseteq$}}
\newcommand{\includes}{\mbox{$\supseteq$}}
\newcommand{\product}{\mbox{$\times$}}
\newcommand{\Power}{\mbox{$\wp$}}
\newcommand{\Powernonempty}{\Power\mbox{$_{+}$}}
\newcommand{\yields}{\mbox{$\mapsto$}}
\newcommand{\comp}{\mbox{$\circ$}}
\newcommand{\emptys}{\mbox{$\emptyset$}}
\newcommand{\emptypres}{\emptys}
\newcommand{\nonabove}{\mbox{$\preceq$}}
\newcommand{\fullynonabove}{\mbox{$\sqsubseteq$}}
\newcommand{\classicallynonabove}{\mbox{$\leq$}}
\newcommand{\carrierof}[1]{\mbox{$|$}#1\mbox{$|$}}
\newcommand{\setof}[1]{\mbox{$\{$}#1\mbox{$\}$}}
\newcommand{\suchthat}{\mbox{$|$}}
\newcommand{\unionall}[3]{\mbox{$\bigcup_{#1\in#2}$}#3}
\newcommand{\widunionall}[3]{\parbox[t]{#3}
                         {\centering\mbox{$\bigcup$}
                          \\ \protect\vspace*{-1.5mm}
                          {\scriptsize #1\of#2}}}
\newcommand{\intsnall}[3]{\mbox{$\bigcap_{#1\in#2}$}#3}
\newcommand{\widintsnall}[3]{\parbox[t]{#3}
                         {\centering\mbox{$\bigcap$}
                          \\ \protect\vspace*{-1.5mm}
                          {\scriptsize #1\of#2}}}
\newcommand{\ifamily}[3]{\setof{#1\mbox{$_{#2}$}}\mbox{$_{#2\in#3}$}}
\newcommand{\indx}{{\it j}}
\newcommand{\anindx}{{\it k}}
\newcommand{\indxset}{{\it J}}
\newcommand{\iso}{\mbox{$\simeq$}}
\newcommand{\lweq}{\mbox{$\leq$}}
\newcommand{\greq}{\mbox{$\geq$}}
\newcommand{\grth}{\mbox{$>$}}
\newcommand{\eqdef}{\mbox{$\stackrel{\it def}{=}$}}
\newcommand{\inverseof}[1]{#1\mbox{$^{-1}$}}
\newcommand{\minverseof}[1]{#1^{-1}}

\newcommand{\seqof}[3]{#1\mbox{$_{1}$}#2\ldots#2#1\mbox{$_{#3}$}}

\newcommand{\seqSubOverof}[3]{#1\raisebox{-2.2mm}
                                     {{\tiny 1}}#2\ldots#2#1\raisebox{-2.2mm}
                                                                 {{\tiny #3}}}

\newcommand{\rpsreq}{(\dag)}
\newcommand{\epsreq}{(\ddag)}
\newcommand{\satinv}{(\S)}
\newcommand{\adequacy}{(\S\S)}
\newcommand{\fulladequacy}{(\S\S\S)}
\newcommand{\finitarity}{(\pounds)}

\newcommand{\upar}{\uparrow}
\newcommand{\Upar}{\Uparrow}
\newcommand{\doar}{\downarrow}
\newcommand{\Doar}{\Downarrow}

\newcommand{\riar}{\rightarrow}
\newcommand{\lear}{\leftarrow}
\newcommand{\Riar}{\Rightarrow}
\newcommand{\Lear}{\Leftarrow}
\newcommand{\leriar}{\leftrightarrow}
\newcommand{\Leriar}{\Leftrightarrow}

\newcommand{\uC}{\underline {\bf C}}
\newcommand{\uD}{\underline {\bf D}}
\newcommand{\uE}{\underline E}
\newcommand{\ugS}{\underline\Sigma}
\newcommand{\uunion}{\underline\cup}
\newcommand{\ubigunion}{\underline\bigcup}
\newcommand{\unE}{\underline E}

\newcommand{\oE}{{\overline E}}
\newcommand{\oP}{{\overline P}}
\newcommand{\oV}{{\overline V}}
\newcommand{\oW}{{\overline W}}
\newcommand{\oU}{{\overline U}}
\newcommand{\oD}{{\overline D}}

\newcommand{\cA}{{\cal A}}
\newcommand{\cB}{{\cal B}}
\newcommand{\cC}{{\cal C}}
\newcommand{\cD}{{\cal D}}
\newcommand{\cE}{{\cal E}}
\newcommand{\cF}{{\cal F}}
\newcommand{\cG}{{\cal G}}
\newcommand{\cH}{{\cal H}}
\newcommand{\cJ}{{\cal J}}
\newcommand{\cL}{{\cal L}}
\newcommand{\cM}{{\cal M}}
\newcommand{\cN}{{\cal N}}
\newcommand{\cO}{{\cal O}}
\newcommand{\cP}{{\cal P}}
\newcommand{\cR}{{\cal R}}
\newcommand{\cS}{{\cal S}}
\newcommand{\cT}{{\cal T}}
\newcommand{\cU}{{\cal U}}
\newcommand{\cV}{{\cal V}}
\newcommand{\cW}{{\cal W}}

\newcommand{\la}{\langle}
\newcommand{\ra}{\rangle}

\newcommand{\sm}{\leadsto}
\newcommand{\sopra}[2]{\stackrel{#1}{#2}}

%definitions for \bf letters
\newcommand{\bb}{{\bf b}}
\newcommand{\bA}{{\bf A}}
\newcommand{\bB}{{\bf B}}
\newcommand{\bC}{{\bf C}}
\newcommand{\bD}{{\bf D}}
\newcommand{\bE}{{\bf E}}
\newcommand{\bF}{{\bf F}}
\newcommand{\bK}{{\bf K}}
\newcommand{\bM}{{\bf M}}
\newcommand{\bN}{{\bf N}}
\newcommand{\bP}{{\bf P}}
\newcommand{\bR}{{\bf R}}
\newcommand{\bS}{{\bf S}}
\newcommand{\bT}{{\bf T}}
\newcommand{\bU}{{\bf U}}
\newcommand{\bX}{{\bf X}}
\newcommand{\bRel}{{\bf Rel}}
\newcommand{\bx}{{\bf x }}

\newcommand{\bCX}{{{\bf C}[X]}}

%definitions for lower case greek letters
\newcommand{\ga}{\alpha}
\newcommand{\gb}{\beta}
\newcommand{\gam}{\gamma}
\newcommand{\gd}{\delta}
\newcommand{\gep}{\varepsilon}
\newcommand{\gz}{\zeta}
\newcommand{\geta}{\eta}
\newcommand{\gth}{\theta}
\newcommand{\gi}{\iota}
\newcommand{\gv}{\nu}
\newcommand{\gk}{\kappa}
\newcommand{\gl}{\lambda}
\newcommand{\gm}{\mu}
\newcommand{\gn}{\nu}
\newcommand{\gx}{\xi}
\newcommand{\gp}{\pi}
\newcommand{\gr}{\rho}
\newcommand{\gs}{\sigma}
\newcommand{\gt}{\tau}
\newcommand{\gu}{\upsilon}
\newcommand{\gph}{\varphi}
\newcommand{\gch}{\chi}
\newcommand{\gps}{\psi}
\newcommand{\go}{\omega}

\newcommand{\gG}{\Gamma}
\newcommand{\gF}{\Phi}
\newcommand{\gD}{\Delta}
\newcommand{\gT}{\Theta}
\newcommand{\gP}{\Pi}
\newcommand{\gX}{\Xi}
\newcommand{\gS}{\Sigma}
\newcommand{\gO}{\Omega}
\newcommand{\gL}{\Lambda}
\newcommand\LT{\mbox{\sf LT}}
\newcommand\sfB{\mbox{\sf B}}
\newcommand\DCAI{\mbox{\sf DCA}_I}
\newcommand\RFA{\mbox{\sf RFA}}
\newcommand\RFAI{\mbox{\sf RFA}_I}
\newcommand\RFAJ{\mbox{\sf RFA}_J}
\newcommand\LAA{\mbox{\sf LAA}}
\newcommand\LAAI{\mbox{\sf LAA}_I}
\newcommand\WLAI{\mbox{\sf LAA}_I^w}
\newcommand\SWLAI{\mbox{\sf LAA}_I^{sw}}
\newcommand\LIEI{\mbox{\sf LIE}_I}
\newcommand\CLAI{\mbox{\sf LAA}_I^c}
\newcommand\CLFAI{\mbox{\sf LFA}_I^c}
\newcommand\SLAI{\mbox{\sf LAA}_I^s}
\newcommand\SSLAI{\mbox{\sf LAA}_I^{cs}}
\newcommand\OLAAI{\mbox{\sf OLAA}_I}
\newcommand\LAAIE{\mbox{\sf LAA}_I^\eta}
\newcommand\LAAJ{\mbox{\sf LAA}_J}
\newcommand\BA{\mbox{\sf BA}}
\newcommand\GSA{\mbox{\sf GSA}}
\newcommand\RGA{\mbox{\sf RGA}}
\newcommand\GA{\mbox{\sf GA}}
\newcommand\K{\mbox{\sf K}}
\newcommand\G{\mbox{\sf G}}

\newcommand\FLA{\mbox{\sf FLA}}
\newcommand\FLAI{\mbox{\sf FLA}_I}
\newcommand\FLAJ{\mbox{\sf FLA}_J}
\newcommand\LFA{\mbox{\sf LFA}}
\newcommand\LFAI{\mbox{\sf LFA}_I}
\newcommand\LFAJ{\mbox{\sf LFA}_J}
\newcommand\CA{\mbox{\sf CA}}
\newcommand\LA{\mbox{\sf LA}}
\newcommand\LM{\mbox{\sf LM}}
\newcommand{\CL}{\mbox{ \sf CL}}

\newcommand\ssk{{\bf k}}
\newcommand\sss{{\bf s}}
\newcommand\ssm{\mbox{\boldmath $m$}}
\newcommand\sse{\mbox{\boldmath $\gep$}}
\newcommand\ssi{{\bf I}}
\newcommand\ssl{\mbox{\boldmath $l$}}
\newcommand\sso{\mbox{\boldmath $1$}}
\newcommand\sst{\mbox{\boldmath $2$}}
\newcommand\ssT{{\mbox{\bf T}}}
\newcommand\ssF{{\mbox{\bf F}}}
\newcommand\ccc{\mbox{\boldmath $c$}}

\newcommand\Con{\mbox{Con}}

\newcommand\sbA{{\mbox{\bf A}}}%
\newcommand\ssbA{{\mbox{\footnotesize\bf A}}}%
\newcommand\sbB{\mbox{\boldmath $B$}}%
\newcommand\sbC{{\mathbf C}}%
\newcommand\sbD{\mbox{\boldmath $D$}}%
\newcommand\sbE{\mbox{\boldmath $E$}}%
\newcommand\sbF{{\mbox{\bf F}}}%
\newcommand\sbG{\mbox{\boldmath $G$}}%
\newcommand\sbH{\mbox{\boldmath $H$}}%
\newcommand\sbI{\mbox{\boldmath $I$}}%
\newcommand\sbJ{\mbox{\boldmath $J$}}%
\newcommand\sbK{{\bold K}}%
\newcommand\sbL{{\bold L}}%
\newcommand\sbM{\mbox{\bf M}}%
\newcommand\sbN{\mbox{\bf N}}%
\newcommand\sbO{{\bold O}}%
\newcommand\sbP{{\bold P}}%
\newcommand\sbQ{{\bold Q}}%
\newcommand\sbR{{\bold R}}%
\newcommand\sbS{{\bold S}}%
\newcommand\sbT{{\bold T}}%
\newcommand\sbU{{\bold U}}%
\newcommand\sbV{{\bold V}}%
\newcommand\sbW{{\bold W}}%
\newcommand\sbX{{\bold X}}%
\newcommand\sbY{{\bold Y}}%
\newcommand\sbZ{{\bold Z}}%

\newcommand\oa{{\overline{a}}}
\newcommand\ob{{\overline{b}}}
\newcommand\oc{{\overline{c}}}
\newcommand\od{{\overline{d}}}
\newcommand\ox{{\overline{x}}}
\newcommand\oy{{\overline{y}}}
\newcommand\ou{{\overline{u}}}
\newcommand\ov{{\overline{v}}}
\newcommand\ow{{\overline{w}}}
\newcommand\oz{{\overline{z}}}
\newcommand\oQ{{\overline{Q}}}
\newcommand\oN{{\overline{N}}}
\newcommand\oM{{\overline{M}}}
\newcommand\oR{{\overline{R}}}
\newcommand\oX{{\overline{X}}}
\newcommand\ogx{{\overline{\gx}}}
\newcommand\oempty{{\overline{\emptyset}}}

\newcommand\twohead{{\twoheadrightarrow}}

\newcommand\up{{\uparrow}}
\newcommand\down{{\downarrow}}
\newcommand\updown{{\updownarrow}}
\newcommand\si{{\leftarrow}}
\newcommand\de{{\rightarrow}}

\newcommand{\bCI}{{\bC[I]}}
\newcommand{\bDI}{{\bD[I]}}

\newcommand\sseq{{\subseteq}}
\newcommand\arrow{{\to}}

\newcommand{\A}{{\mathfrak{A}}}
\newcommand{\di}{{\diamondsuit}}

\newcommand{\Th}{{\mbox{Eq}}}
\newcommand{\Or}{{\mbox{Ord}}}

\newcommand{\sq}{{\sqsubseteq}}

\newcommand{\BS}{{\cal B}_{Sc}}

\maketitle
\thispagestyle{empty}
\vspace*{-2.7mm}
\begin{abstract} 
A longstanding open problem 
is whether there exists a non-syntactical model of the untyped
$\lambda$-calculus whose theory is exactly the least $\lambda$-theory  $\lambda\beta$.
In this paper we investigate the more general question of whether 
the equational/order theory of a model of the untyped  $\lambda$-calculus can be recursively enumerable (r.e. for brevity).
We introduce a notion of \emph{effective model} of  $\lambda$-calculus, which covers in particular all the models
individually introduced  in the literature.
 We prove that
the order theory of an effective model is never r.e.; from
this it follows that its equational theory cannot be $\lambda\beta$, $\lambda\beta\eta$.
We then show that 
no effective model living in the stable or strongly stable semantics 
has an r.e. equational theory.
Concerning Scott's semantics, we investigate the class of graph models and prove
that  no order theory of a graph model can be r.e., and that
there exists an effective graph model whose  equational/order  theory is the minimum one.
Finally, we show that the class of graph models enjoys a kind of downwards L\"owenheim-Skolem theorem.
\\

\end{abstract}

\vspace*{-6.0mm}
\section{Introduction}
\vspace*{-1.0mm}
Lambda theories are equational extensions of the untyped $\lambda$-calculus closed under derivation.
They arise by syntactical or semantic considerations.
Indeed, a  $\gl$-theory may correspond to a possible operational (observational) semantics of $\lambda$-calculus, as well as
it  may be induced by a model of $\lambda$-calculus through the kernel congruence relation of the interpretation function. 
Although researchers have mainly focused their interest on a limited number of them, the class of $\gl$-theories
constitutes a very rich and complex structure (see \cite{Bare,Berli,Berli2}). 

Topology is at the center of the known approaches to giving models of the untyped $\lambda$-calculus. 
After the first model, found by Scott in 1969 in the category of complete lattices and 
Scott continuous functions, a large number of mathematical models for $\lambda$-calculus,
arising from syntax-free constructions, have been introduced in various 
categories of domains  and were classified into semantics according to 
the nature of their representable functions, see e.g. \cite{Bare,Berli,Plot}. 
Scott continuous semantics \cite{Sc} is given in the category whose objects 
are complete partial orders and morphisms are Scott continuous functions.
The stable semantics (Berry \cite{Ber}) and the strongly stable 
semantics (Bucciarelli-Ehrhard \cite{BuEh})
are refinements of the continuous semantics, introduced to 
capture the notion of ``sequential" Scott continuous function. 
In each of these semantics it is possible to build up $2^{\aleph_0}$  models inducing pairwise distinct  $\gl$-theories
\cite{Ke,Ker}.
Nevertheless, 
all are equationally
{\em  incomplete} (see \cite{HR,BG,Salics,Sacm}) in the sense that they do not represent all possible 
consistent $\lambda$-theories.
It is interesting to note that there are very few known equational theories of $\gl$-models living in these semantics
that can be described syntactically:
namely, the theory of B\"ohm trees and variants of it. None of these theories is r.e.

Berline  has raised in \cite{Berli}
the natural question of whether, given a class  of models of $\lambda$-calculus, 
there is a minimum  $\gl$-theory represented by it.
This question relates to the longstanding open problem proposed by Barendregt 
about the existence of a continuous model or, more generally, of a non-syntactical model of 
$\gl\gb$ ($\gl\gb\eta$). 
Di Gianantonio, Honsell and Plotkin \cite{DHP} have shown that
Scott continuous semantics admits a minimum theory, at least if we restrict to extensional models. 
Another result of  \cite{DHP}, in the same spirit, is the 
construction of an extensional model whose theory is $\lambda\beta\eta$, a
fortiori minimal, in a weakly-continuous semantics.
However, the construction of this model starts from the term model of $\gl\gb\eta$, and hence
it cannot be seen as having a purely non syntactical presentation. 
More recently, Bucciarelli and Salibra \cite{BSa,BSa2}
have shown that the class of graph models admits a minimum  $\gl$-theory different from $\gl\gb$. 
Graph models, isolated in the seventies by Plotkin,
Scott and Engeler (see e.g. \cite{Bare}) within the continuous semantics, have proved 
useful for showing the consistency of extensions
of $\lambda$-calculus and for studying operational features of
$\lambda$-calculus (see \cite{Berli}).

In this paper we  investigate the related question of whether 
the equational theory of a model  
can be recursively enumerable (r.e. for brevity).
As far as we know, this problem was first raised in
\cite{Berli2}, where it is conjectured that no graph model can have an r.e. theory.
But we expect that this could indeed be true for all models living in the continuous semantics, 
and its refinements.

We find it natural to concentrate 
on models with built-in effectivity properties.
It seems indeed reasonable to think that, if effective models do not even succeed to have an r.e. theory, 
then the other ones have no chance to succeed.
Another justification for considering effective models comes from a previous result obtained
for typed $\gl$-calculus. Indeed, it was proved in \cite{BB} that
there exists a non-syntactical model of Girard's system $F$ whose theory is $\gl\gb\eta$.
 This model lives in Scott's continuous semantics, and can  easily be
 checked to be ``effective'' in the same spirit as in the present paper (see \cite[Appendix C]{BB}
for a sketchy presentation of the model).

Starting from the known notion of an effective domain, 
we introduce a general notion of an
{\em effective model} of $\lambda$-calculus 
and we study the main properties of these models.
Effective models are omni-present in the continuous, stable and strongly stable semantics.
In particular, all the models which have been introduced individually in the literature can easily be proved 
effective\footnote{As far as we know, only Giannini and Longo \cite{GL} have introduced a  notion of
an effective model; but their definition is
\emph{ad hoc} for two particular models
(Scott's $P_\go$ and Plotkin's $T_\go$) and their results depend on the fact that
these models have a very special (and well known) common theory.}.

%%%%%%%%%%%%%%%%%%%%%%%%%%%%%%%%%%%%
The following are the main results of the paper:

\begin{enumerate}
\item Let $\cD$ be an effective model of $\gl$-calculus. Then:
\begin{itemize}
\item[(i)] The order theory $\Or(\cD)$ of $\cD$ is not r.e.
\item[(ii)] The equational theory $\Th(\cD)$ of $\cD$ is not the theory $\gl\gb$ ($\gl\gb\eta$).
\item[(iii)] If for some $\gl$-term $M$ there are only finitely many  $\gl$-definable elements below the interpretation of $M$
(e.g. if $\bot \in \cD$ is $\gl$-definable), then $\Th(\cD)$ is not r.e.
\end{itemize}
\end{enumerate}
Concerning the existence of a non-syntactical effective model with an r.e. equational theory, we are able to give 
a definite negative answer for all (effective) stable and strongly stable
models: 

\begin{itemize}
\item[2.] No effective model living in the stable or strongly 
stable semantics has an r.e. equational theory.
\end{itemize}
Concerning Scott continuous semantics, the problem looks much more difficult. We concentrate here
on the class of graph models (see \cite{Berli2,BeSa,BSa,BSa2,BSa3} for earlier investigation of this class)
and show the following results:

\begin{itemize}
\item[3.] Let $\cD$ be an arbitrary graph model. Then:
\begin{itemize}
\item[(i)] The order theory $\Or(\cD)$ of $\cD$ is not r.e.
\item[(ii)] If $\cD$ is freely generated  by a finite ``partial model'',
then the equational theory $\Th(\cD)$ of $\cD$ is not r.e.
\end{itemize}
\end{itemize}

\begin{itemize}
\item[4.] There exists an effective graph model whose
equational/order theory is minimal among all theories of graph models.
\end{itemize}

\begin{itemize}
\item[5.] (L\"owenheim-Skolem theorem for graph models) 
Every equational/order graph theory (where ``graph theory'' means ``theory of a graph model'') 
is the theory of a graph model having a carrier set of minimal cardinality.
\end{itemize}
The last result positively 
 answers Question 3 in \cite[Section 6.3]{Berli} for the class of graph models. 
%and has the consequence that every graph theory (we know from Kerth \cite{Ke} that there exists a continuum of them)
%is the theory of a graph model whose carrier set is the set of natural numbers.

The central technical device used in this paper is Visser's result \cite{Vi} stating that 
the  complements of $\gb$-closed r.e. sets of $\gl$-terms enjoy the finite intersection property  
(see Theorem~\ref{inter-prop}).

%%%%%%%%%%%%%%%%%%%%%%%%%%%%%%%%%%%%

\vspace*{-2.0mm}
\section{Preliminaries}\label{pre}
\vspace*{-1.0mm}
To keep this article self-contained, we summarize some definitions and 
results concerning  $\gl$-calculus that we need in the subsequent part of the paper.
With regard to the lambda calculus we follow the notation and terminology of \cite{Bare}.

We denote by $\nat$ the set of natural numbers.
A set $A\subseteq\nat$ is \emph{recursively enumerable} (r.e. for short) if it is the domain of a 
partial recursive function. 
The complement of a recursively enumerable set  is called a \emph{co-r.e.} set.
If both $A$ and its complement are r.e., $A$ is called \emph{decidable}.
We will denote by $\RE$ the collection of all r.e. subsets of $\nat$.

A numeration of a set $A$ is a map from $\nat$ onto $A$.
$\cW: \nat\to\RE$ denotes   the usual numeration of r.e. sets
(i.e., $\cW_n$ is the domain of the $n$-th computable function $\phi_n$).

\vspace*{-2.0mm}
\subsection{Lambda calculus and lambda models}\label{lambdacalculus}
\vspace*{-1.0mm}
$\gL$ and $\gL^o$ are, respectively, the set of $\gl$-terms and of closed $\gl$-terms. 
Concerning specific $\gl$-terms we set: 
$$\ssi \equiv \gl x.x;\quad \ssT \equiv \gl xy.x;\quad \ssF\equiv \gl xy.y;\quad \gO \equiv (\gl x.xx)(\gl x.xx).$$

A set $X$ of $\gl$-terms is {\em trivial} if either $X = \emptyset$ or $X= \gL$.

We denote $\ga\gb$-conversion by $\gl\gb$. 
A {\em  $\gl$-theory} $\cT$ is a congruence on $\gL$ (with respect to the operators of abstraction
and application) which contains $\gl\gb$.
We write $M =_\cT N$ for $(M, N) \in \cT$.
If $\cT$ is a $\gl$-theory, then $[M]_\cT$ denotes the set $\{ N: N=_\cT M\}$.
A $\gl$-theory $\cT$ is: \emph{consistent} if $\cT \neq \gL\times\gL$;
\emph{extensional} if it contains  the equation $\ssi = \gl xy.xy$;
 \emph{recursively enumerable} if the set of G\"odel 
numbers of all pairs of $\cT$-equivalent $\gl$-terms is r.e.
Finally, $\gl\gb\eta$ is the least extensional $\gl$-theory.

Solvable $\gl$-terms can be characterized as follows:
a $\gl$-term $M$ is solvable if, and only if, it has a {\em head normal form},
that is, $M =_{\gl\gb} \gl x_{1} \ldots x_{n}.yM_1\dots M_k$
for some $n,k \geq 0$ and $\gl$-terms $M_1,\dots, M_k$.
$M\in\gL$ is {\em unsolvable} if it is not solvable. 

The $\gl$-theory $\cH$, generated by equating all the unsolvable $\gl$-terms, is consistent 
by \cite[Theorem~16.1.3]{Bare}. 
A $\gl$-theory $\cT$ is \emph{sensible} if $\cH\subseteq\cT$, while it is
 \emph{semi-sensible} if it contains no equations of the form $U = S$ where $S$ 
is solvable and $U$ unsolvable. 
Consistent sensible theories are semi-sensible (see \cite[Cor.~4.1.9]{Bare})
and are never r.e. (see \cite[Section~17.1]{Bare}).

It is well known \cite[Chapter~5]{Bare} that a model of $\lambda$-calculus ({\em $\gl$-model}, for short)  can be defined
as a reflexive object in a ccc (Cartesian closed category) 
$\bC$, 
that is to say a triple $(D,\App,\Abs)$ such that $D$
is an object of $\bC$ and 
$\App: D\to[ D\to D]$, $\Abs:[ D\to D]\to D$ are morphisms such that $\App\comp\Abs=id_{[D\to D]}$. 
In the following we will mainly be interested in Scott's ccc of
cpos and Scott continuous functions (\emph{continuous semantics}), but
we will also draw conclusions for Berry's ccc of $DI$--domains and stable functions
(\emph{stable semantics}), and for Ehrhard's ccc of $DI$-domains
with coherence and strongly stable functions between them (\emph{strongly stable semantics}). 
We recall that $DI$-domains are special Scott domains, and that Scott domains are special cpos
  (see, e.g.,  \cite{SLG}).

Let $D$ be a cpo. The partial order of $D$ will be denoted by $\sq_D$. 
We let $Env_D$ be the set of environments $\gr$ mapping the set $Var$ of variables of  $\gl$-calculus into $D$.
For every $x\in Var$ and $d\in D$ we denote by $\gr[x:=d]$ the environment $\gr'$ which coincides with $\gr$, except on $x$,
where $\gr'$ takes the value $d$.
A reflexive cpo $D$ generates a $\gl$-model $\cD = (D,\App,\Abs)$  with the interpretation 
of a $\gl$-term defined as follows:
$$x^\cD_\gr= \gr(x);\ (MN)^\cD_\gr = \App(M^\cD_\gr)(N^\cD_\gr);\ 
(\gl x.M)^\cD_\gr = \Abs(f),$$
where $f$ is defined by $f(d)=M^\cD_{\gr[x:=d]}$ for all $d\in D$.
We write $M^\cD$ for $M^\cD_\gr$ if $M$ is a closed $\gl$-term.
In the following $\App(d)(e)$ will also be written $d\cdot e$ or $de$.

Each $\gl$-model $\cD$ induces a  $\gl$-theory, denoted here by $\Th(\cD)$,
and called {\em the equational theory of $\cD$}. Thus,
$M = N \in \Th(\cD)$ if, and only if, $M$ and $N$ have the same interpretation in $\cD$.
A reflexive cpo $\cD$ induces also an {\em order theory} 
$\Or(\cD)=\{ M\ \sq\ N : M^\cD_\gr\ \sq_D\ N^\cD_\gr\ \mbox{for all environments $\gr$}\}$.

\vspace*{-2.0mm}
\subsection{Effective domains}\label{domain}
\vspace*{-1.0mm}
A triple $\cD= (D,\sq_D, d)$ is called an \emph{effective domain} if $(D,\sq_D)$ is a Scott domain and
$d$ is a numeration
of the set $K(\cD)$ of its compact elements 
such that the relations ``$d_m$ and $d_n$ have an upper bound'' and
``$d_n = d_m \sqcup d_k$" are both decidable (see, e.g.,  \cite[Chapter~10]{SLG}).

We recall that  an element $v$ of an effective domain $\cD$ is said \emph{r.e.} (\emph{decidable})
 if the set $\{ n : d_n \sq_D v\}$ is r.e. (decidable); 
we will write $\cD^{r.e.}$ ($\cD^{dec}$) for the set of r.e. (decidable) elements of $\cD$.
The set $K(\cD)$ of compact elements is 
included within $\cD^{dec}$.
Using standard techniques of recursion theory it is possible to get in a uniform way a numeration 
$\xi: \nat\to \cD^{r.e.}$ which is \emph{adequate} in the sense that the relation $d_k\ \sq_D\ \xi_n$ is r.e. in $(k,n)$ 
and the inclusion mapping $\iota:K(\cD)\to \cD^{r.e.}$ is computable w.r.t. $d,\xi$.

The full subcategory $\bold{ED}$ of the category of Scott-domains 
with effective domains as objects and continuous functions as morphisms is a ccc.

A continuous function $f:D\to D'$ is an r.e. element in the effective 
domain of Scott continuous functions (i.e., $f\in [\cD\to \cD']^{r.e.}$)  if, and only if, its restriction 
$f\hspace{-4pt}\restr:\cD^{r.e.}\to \cD'^{r.e.}$ 
is {\em computable w.r.t. $\xi,\xi'$}, i.e., there is a computable map $g:\nat\to\nat$ such that
$f(\xi_n) = \xi'_{g(n)}$. In such a case we say that $g$ \emph{tracks} $f$.

\vspace*{-2.0mm}
\subsection{Graph models}\label{graph-models}
\vspace*{-1.0mm}
The class of graph models belongs to Scott continuous semantics (see \cite{Berli2} for a complete survey on this class of models).
Historically, the first graph model was Scott's $P_\go$, which is also
known in the literature as ``the graph model". ``Graph'' referred to the fact
that the continuous functions were encoded in the model via (a sufficient fragment of)
their graph.

As a matter of notation, for every set $G$,
$G^*$ is the set of all finite subsets of 
$G$, while $\cP(G)$ is the powerset of $G$. 

\begin{definition} 
A {\em graph model} $\cG$ is a pair $(G,c_\cG)$, where $G$ is an infinite set, called the {\em carrier set} of $\cG$, and 
$c_\cG: G^*\times G \to G$ is an injective total function.
\end{definition}

Such pair $\cG$ generates the reflexive cpo $(\cP(G),\subseteq, \Abs, \cF)$,
where $\Abs$ and $\App$ are defined as follows, for all 
$f\in [\cP(G)\to \cP(G)]$ and  $X,Y\subseteq G$:
$\Abs(f) = \{ c_\cG(a, \ga) : \ga \in f(a)\ \mbox{and $a\in G^*$}\}$ and
$\App(X)(Y) = \{ \ga\in G : (\exists a\subseteq Y)\ c_\cG(a,\ga)\in X\}$.
For more details we refer the reader to Berline \cite{Berli}.

The interpretation of a $\gl$-term $M$ into a $\gl$-model
has been defined in Section~\ref{lambdacalculus}. However, in this context we can make explicit 
the interpretation $M^\cG_\gr$ of a $\gl$-term $M$ as follows:
$$
(MN)^\cG_\gr= \{ \ga: 
(\exists a\subseteq N^\cG_\gr)\ c_\cG(a,\ga)\in M^\cG_\gr\}; 
(\gl x.M)^\cG_\gr = \{ c_\cG(a,\ga) :  \ga\in M^\cG_{\gr[x:=a]}\}.$$

We turn now to the interpretation of $\gO$ in graph models (the details of the proof
are, for example, worked out in \cite[Lemma 4]{BeSa}). 

\begin{lemma}\label{omegan} $\ga\in\gO^\cG$ if, and only if, there is $a\subseteq  (\gl x.xx)^\cG$ such that 
 $c_\cG(a,\ga)\in a$.
\end{lemma}

In the following we use the terminology  ``{\em graph theory}'' as a shorthand for ``theory of a graph model''.
It is well known that the equational graph theories are never extensional and that
there exists a continuum of them (see \cite{Ke}). 
In \cite{BSa,BSa2} the existence of a minimum equational graph theory was proved and it was also
shown that this minimum theory is different from $\gl\gb$.

The completion method for building graph models from ``partial pairs'' was initiated
by Longo in \cite{Lon} and developed on a wide scale by Kerth in 
\cite{Ke,Ker}.

\begin{definition} 
A {\em partial pair} $\cA$ is given by a set $A$ and
by a partial, injective function $c_\cA:A^*\times A\arrow A$.
\end{definition}

A partial pair is {\em finite} if $A$ is finite, and is a graph model if $c_\cA$ is total. 

The interpretation of a $\gl$-term  in a partial pair $\cA$ is defined in the obvious way: 
$(MN)_{\gr }^{\cA}=\{\,\ga \in A\,:\,(\exists a\subseteq N_{\gr
}^{\cA})\ [(a,\ga )\in dom(c_\cA)\wedge c_\cA(a,\ga )\in M_{\gr }^{\cA}]\}$; 
$(\gl x.M)_{\gr }^{\cA}=\{\,c_\cA(a,\ga )\in \cA\,:\,(a,\ga )\in
dom(c_\cA)\wedge \ga \in M_{\gr \lbrack x:=a]}^{\cA}\,\}$.

\begin{definition}\label{completion}
Let $\cA$ be a partial pair. The {\em completion} of $\cA$ 
is the graph model $\cE_\cA = (E_\cA, c_{\cE_\cA})$  defined as follows:

\begin{itemize}
\item $\ E_\cA=\bigcup_{n\in\nat} E_n$, where $E_0=A$ and $E_{n+1}=
E_n\ \cup\ ((E_n^*\times E_n)- dom(c_\cA))$.

\item $\ $Given $a\in E_\cA^*$, $\ga\in E_\cA$,
\vspace*{-2.0mm}
$$c_{\cE_\cA}(a,\ga)=\left\{\begin{array}{ll}
c_\cA(a,\ga)&\mbox{if}\ c_\cA(a,\ga)\ \mbox{is defined} \\
(a,\ga)& \mbox{otherwise}
\end{array}
\right.
$$
\end{itemize}
\end{definition}

A notion of {\em rank} can be naturally defined on the  completion $\cE_\cA$ of a partial pair $\cA$.
The elements of $A$ are the elements of rank $0$, while an element $\ga \in E_\cA-A$
has rank $n$ if $\ga \in E_n$ and $\ga\not\in E_{n-1}$.

Let $\cA$ and $\cB$ be two partial pairs. A {\em morphism} from $\cA$ into $\cB$
is a map $f: A\to B$ such that 
$(a,\ga) \in dom(c_\cA)$ implies $(fa,f\ga)\in dom(c_\cB)$ and, in such a case $f(c_\cA(a,\ga)) = c_\cB(fa,f\ga)$.
Isomorphisms and automorphisms can be defined in the obvious way.
$Aut(\cA)$ denotes the group of automorphisms of the partial pair $\cA$.

\begin{lemma}\label{mor} Let $\cG,\cG'$ be graph models and $f:\cG\to\cG'$ be a morphism. 
If $M\in\gL$ and $\ga\in M^{\cG}_\gr$, then $f\ga \in M^{\cG'}_{f\circ \gr}$.
\end{lemma}

%%%%%%%%%%%%%%%%%%%%%%%%%%%%%%%%%%%%%%%%%%%%%%%%%%

\vspace*{-1.0mm}
\subsection{Co-r.e. sets of lambda terms }\label{sec:Vissertop}
%\vspace*{-1.0mm}
In this section we recall the main properties of recursion theory concerning  $\gl$-calculus that
will be applied in the following sections.

An r.e. (co-r.e.) set  of $\gl$-terms closed under $\gb$-conversion will be called a {\em $\gb$-r.e.} 
({\em $\gb$-co-r.e.}) set.

The following theorem is due to Scott (see \cite[Thm. 6.6.2]{Bare}).

\begin{theorem}\label{both} A set of $\gl$-terms which is both $\gb$-r.e. and $\gb$-co-r.e. is trivial.
\end{theorem}

\begin{definition} A family $X = (X_i : i\in I)$ of sets has the {\em FIP} (finite intersection property) if 
$X_{i_1}\cap\dots\cap X_{i_n}\neq \emptyset$ for all $i_1,\dots,i_n\in I$.
\end{definition}

Visser (see \cite[Ch. 17]{Bare} and \cite[Thm.~2.5]{Vi}) has shown  
that the topology on $\gL$ generated by the $\gb$-co-r.e. sets of $\gl$-terms is hyperconnected (i.e., the intersection
of two non-empty open sets is non-empty). In other words:

\begin{theorem}\label{inter-prop} The family of all non-empty $\gb$-co-r.e. subsets of $\gL$ has the FIP.
\end{theorem}

\begin{corollary}\label{uns} Every non-empty $\gb$-co-r.e. set of $\gl$-terms contains a non-empty
$\gb$-co-r.e. set of unsolvable $\gl$-terms.
\end{corollary}

\begin{proof} The set of all unsolvable $\gl$-terms is $\gb$-co-r.e.
The conclusion follows from Theorem~\ref{inter-prop}.
\end{proof}

%%%%%%%%%%%%%%%%%%%%%%%%%%%%%%%%%%%%%%%%%%%%%%%%%%

\vspace*{-2.0mm}
\section{Effective lambda models}\label{section-effective}
\vspace*{-1.0mm}
In this section we introduce the notion of an effective $\gl$-model and we study the main properties of these models.
We show that the order theory of  an effective $\gl$-model  is not r.e. and that its equational theory  is
different from $\gl\gb, \gl\gb\eta$.
Effective $\gl$-models are omni-present in the continuous, stable and strongly stable semantics 
(see Section~\ref{model-examples}). 
In particular, all the $\gl$-models which have been introduced individually in the literature, to begin with
Scott's $\mathcal{D}_{\infty}$, 
can  easily be  proved effective. 

The following natural definition is enough to force
the interpretation function of $\gl$-terms to be computable from $\gL^o$ into $\cD^{r.e.}$. 
However, other results of this paper will need a more powerful notion. That is the reason why
we only speak of ``weak effectivity'' here.

\begin{definition}\label{weakly-effective} A $\gl$-model is called \emph{weakly effective}
if it is a reflexive object $(\cD,\App,\Abs)$ in the category $\bold{ED}$ and,
 $\App\in [\cD\to[\cD\to \cD]]$ and $\Abs\in [[\cD\to \cD]\to \cD]$ are r.e. elements.
\end{definition}

In the following a weakly effective $\gl$-model $(\cD,\App,\Abs)$ will be denoted
by $\cD$.

We fix bijective effective numerations $\nu_{\gL}:\nat\to\gL$ of the set of $\gl$-terms and
 $\nu_{var}:\nat\to Var$ of the set  of variables of $\lambda$-calculus.
In particular this gives to the set $Env_D$ of all environments a structure of effective domain. 
$\gL_\bot = \gL\cup\{\bot\}$ is the usual flat domain of $\gl$-terms. 
The element $\bot$ is always interpreted as $\bot_D$
in a cpo $(D,\sq_D)$.

\begin{proposition}\label{sketch} Let $\cD$ be a weakly effective $\gl$-model.  
Then the function $f$ mapping $(\gr,M)\mapsto M^\cD_\gr$ is an element of $[Env_D\times\gL_\bot \to \cD]^{r.e.}$.
\end{proposition}

\begin{proof}\bold{ (Sketch)} By structural induction on $M$ it is possible to
show the existence of a partial computable map tracking $f$. 
The only difficult case is $M\equiv \gl x.N$. 
Since $\Abs$ is r.e. it is sufficient to prove that the function $g:e\mapsto N_{\gr[x:=e]}^\cD$ is also r.e. 
Once shown that $h:(\gr,x,e)\mapsto \gr[x:=e]$ is r.e.,
from the induction hypothesis  it follows that the function 
$g'(\gr,x,e) = f(h(\gr,x,e),N)$ is r.e. Then by applying the s-m-n theorem of recursion theory
 to the computable function tracking $g'$ we obtain a computable function 
tracking $g$, which is then r.e.
\end{proof}

\begin{notation} We define for any $e\in D$ and $M\in \gL^o$: 

(i) $e^- \equiv \{ P\in\gL^o : P^\cD\ \sq_D\ e\}$; 

(ii) $M^- \equiv \{ P\in\gL^o : P^\cD\ \sq_D\ M^\cD\}$.
\end{notation}

\begin{corollary}\label{nino}  If $e\in \cD^{dec}$, then $e^-$ is a $\gb$-co-r.e. set of $\gl$-terms.
\end{corollary}

\begin{proof} Let $\gr\in (Env_D)^{r.e.}$ be an environment. 
By Proposition~\ref{sketch} there is a computable map $\phi$ tracking 
the interpretation function $M\mapsto M^\cD_\gr$ of $\gl$-terms from $\gL$ into $\cD^{r.e.}$
with respect to the effective numeration $\nu_{\gL}$ of $\gL$ and an adequate numeration $\gx$
of  $\cD^{r.e.}$.
From $e\in \cD^{dec}$ it follows that the set $X = \{ n : \gx_n\ \sq_\cD\ e\}$ is co-r.e.
This implies that the set $\phi^{-1}(X)$, which is the set of the codes of the elements of
 $\{ M\in \gL : M^\cD_\gr\ \sq_\cD\ e\}$,
is also co-r.e. We get the conclusion because $\gL^o$ is a decidable subset of $\gL$.
\end{proof}

\begin{definition}\label{effective}
A weakly effective $\gl$-model $\cD$ is called {\em effective} if it satisfies the following two further conditions:
\begin{itemize}
\item [(i)] If $d\in K(\cD)$ and $e_i\in \cD^{dec}$, then $de_1\dots e_n\in \cD^{dec}$.
\item [(ii)] If $f\in [\cD\to \cD]^{r.e.}$ and $f(e)\in \cD^{dec}$ for all %compact elements 
$e\in K(\cD)$, then $\Abs(f) \in \cD^{dec}$.
\end{itemize}
\end{definition}

An environment $\gr$ is {\em compact} in the effective domain $Env_D$ (i.e., $\gr\in K(Env_D)$)
 if $\gr(x) \in K(\cD)$ for all variables $x$
and $\{ x: \gr(x) \neq \bot_D\}$ is finite.

\begin{notation}\label{not} We define:
$\gL^{dec}_\cD \equiv\{ M\in\gL : \mbox{$M^\cD_\gr\in \cD^{dec}$ for all $\gr\in K(Env_D)\}$}.$
\end{notation}

\begin{theorem}\label{thm} Suppose $\cD$ is an effective $\gl$-model.
Then the set  $\gL^{dec}_\cD$ is closed under the following rules:
\begin{enumerate}
\item $x\in \gL^{dec}_\cD$ for every variable $x$.
\item $M_1,\dots, M_k\in \gL^{dec}_\cD\ \Rightarrow\ yM_1\dots M_k \in \gL^{dec}_\cD$.
\item  $M\in \gL^{dec}_\cD\ \Rightarrow\ \gl x.M \in \gL^{dec}_\cD$.
\end{enumerate}
In particular, $\gL^{dec}_\cD$ contains all the $\gb$-normal forms.
\end{theorem}

\begin{proof} Let $\gr \in K(Env_D)$. We have three cases.\\
(1) $x^\cD_\gr = \gr(x)$ is compact, hence it is decidable.\\
(2) %Suppose $M\equiv yM_1\dots M_k$ with $M_i\in \gL^{dec}_\cD$ for all $i$. 
By definition $(yM_1\dots M_k)^\cD_\gr = \gr(y)(M_1)^\cD_{\gr}\dots (M_k)^\cD_{\gr}$.
Hence the result follows from Definition~\ref{effective}(i),
  $\gr(y)\in K(\cD)$ and $(M_i)^\cD_{\rho}\in \cD^{dec}$.\\
(3) By definition we have that $(\gl x.M)^\cD_\gr = \Abs(f)$,
where $f(e)= M^\cD_{\gr[x:= e]}$ for all $e\in D$. 
Note that $\gr[x:=e]$ is also compact for all $e\in K(D)$.
Hence the conclusion follows from $M^\cD_{\gr[x:= e]}\in D^{dec}$ ($e\in K(\cD)$), 
Definition~\ref{effective}(ii) and $f\in [\cD\to \cD]^{r.e.}$.
\end{proof}

Recall that $\Th(\cD)$ and $\Or(\cD)$ are respectively the equational theory and the order theory of $\cD$.

\begin{theorem}\label{cor} 
Let $\cD$ be an effective $\gl$-model,
and let $M_1, \dots M_k\in \gL^{dec}_\cD$ ($k\geq 1$) be closed terms. Then we have:
\begin{itemize}
\item[(i)] $M_1^-\cap\dots\cap M_k^-$  is a $\gb$-co-r.e. 
set, which contains a non-empty $\gb$-co-r.e. set of unsolvable terms.
\item[(ii)] If $e\in \cD^{dec}$ and $e^-$ is non-empty and finite modulo $\Th(\cD)$, then $\Th(\cD)$
is not r.e. (in particular, if $\bot_\cD^-\neq\emptyset$ then $\Th(\cD)$
is not r.e.).
\item[(iii)] $\Or(\cD)$ is not r.e.
\item [(iv)] $\Th(\cD) \neq \gl\gb, \gl\gb\geta$.
\end{itemize}
\end{theorem}

\begin{proof} (i) By Theorem~\ref{thm}, Corollary~\ref{uns}, Corollary~\ref{nino} and the FIP.

(ii) By Corollary~\ref{nino} we have that  $e^-$ is a $\gb$-co-r.e. set of closed $\gl$-terms.
The conclusion follows because $e^-$ is non-empty and finite modulo $\Th(\cD)$.

(iii) Let $M\in\gL^{dec}_\cD$ be a closed term. If $\Or(\cD)$ were r.e., then
we could enumerate the set $M^-$. However, by (i) 
this set is non-empty and $\gb$-co-r.e. By Theorem~\ref{both} it follows that $M^- = \gL^o$.
By the arbitrariness of $M$, it follows that
$\ssT^- = \ssF^-$. Since $\ssF \in \ssT^-$ and
conversely we get $\ssF=\ssT$ in $\cD$, contradiction.

(iv) Because of (iii), if $\Th(\cD)$ is r.e. then $\Or(\cD)$
strictly contains $\Th(\cD)$. Hence the conclusion follows from Selinger's
result stating that in any partially ordered $\gl$-model, whose theory is 
$\gl\gb$, the interpretations of distinct closed terms are incomparable \cite[Corollary~4]{Sel}.
Similarly for $\gl\gb\eta$.
\end{proof}

%%%%%%%%%%%%%%%%%%%%%%%%%%%%%%%%%%%%%%%%%%%%

\vspace*{-2.0mm}
\section{Can effective $\gl$-models have an r.e. theory?}\label{model-examples}
\vspace*{-1.0mm}
In this section we give a sufficient condition for a wide class of graph models to be effective and show that 
no effective graph model generated freely by a partial pair, which is finite 
modulo its group of automorphisms, can have an r.e. equational theory. 
Finally, we show that no effective $\gl$-model living in the stable or strongly stable semantics 
can have an r.e. equational theory.

In Section~\ref{skolem} we will show that every equational/order graph theory is the theory of a graph model $\cG$
whose carrier set is the set $\nat$ of natural numbers. 
In the next theorem we characterize the effectivity of these models. 

\begin{theorem}\label{effective-graph} Let $\cG$ be a graph model  such that, after encoding, $G=\nat$
and $c_\cG$ is a computable map.
Then $\cG$ is weakly effective. Moreover, $\cG$ is effective under the further hypothesis that $c_\cG$ has a decidable range.
\end{theorem}

\begin{proof} 
It is easy to check, using the definitions given in Section \ref{graph-models}, 
that $\App,\Abs$ are r.e. in their respective domains and that condition (i) of Definition~\ref{effective} is  satisfied.
Then $\cG$ is weakly effective.
Moreover, Definition~\ref{effective}(ii) holds under the hypothesis that the range of $c_\cG$ is decidable. 
\end{proof}

Completions of partial pairs have been extensively studied in the literature. They are useful for solving
equational and inequational constraints (see \cite{Berli,Berli2,BSa2,BSa3}).
In \cite{BSa3} Bucciarelli and Salibra have recently proved that the theory of the completion of a 
 partial pair which is not a graph model is semi-sensible.
The following theorem shows, in particular, that the theory of the completion of a finite partial pair is not r.e.

\begin{theorem}\label{effective-partial} 
Let $\cA$ be a partial pair such that $A$ is finite or  equal to $\nat$ after encoding,
 and $c_\cA$ is a computable map with a decidable domain. Then we have:
\begin{itemize}
\item[(i)]  The completion $\cE_\cA$ of $\cA$ is weakly effective;
\item[(ii)] If the range of $c_\cA$ is decidable, then $\cE_\cA$ is effective;
\item[(iii)] If $\cA$ is finite modulo its group of automorphisms (in particular, if $A$ is finite), 
then $\Th(\cE_\cA)$ is not r.e.
\end{itemize}
\end{theorem}

\begin{proof} Since $A$ is finite or equal to  $\nat$
we have that $E_\cA$ is also decidable (see Definition~\ref{completion}). Moreover, the map
$c_{\cE_\cA}: E_\cA^*\times E_\cA\to E_\cA$ is computable, because it is an extension of a computable function $c_\cA$
with decidable domain, and it is the identity on the decidable set $(E_\cA^*\times E_\cA)- dom(c_\cA)$.
Then (i)-(ii) follow from Theorem~\ref{effective-graph}.

Clearly $A$ is a decidable subset of $E_\cA$; then by Corollary~\ref{nino} 
the set $A^-$ is a $\gb$-co-r.e. set of $\gl$-terms. We now show that this set is non-empty because
$\gO^{\cE_\cA}\subseteq A$.
By Lemma~\ref{omegan} we have that $\ga\in\gO^{\cE_\cA}$ implies that $c_{\cE_\cA}(a,\ga)\in a$ for some
$a\in E_\cA^*$. Immediate considerations on the rank show 
that this is only possible if $(a,\ga)\in dom(c_\cA)$, which forces $\ga\in A$.

The orbit of $\ga\in A$ modulo $Aut(\cA)$ is defined by
 $O(\ga) = \{\gth(\ga) : \gth\in Aut(\cA)\}$.

We now show that, if the set of orbits of $\cA$ has cardinality $k$ for some $k\in\nat$, then 
the cardinality of $A^-$ modulo $\Th(\cE_\cA)$ is less than or equal to $2^k$. 
Assume $p\in M^{\cE_\cA}\subseteq A$. Then by Lemma~\ref{mor}
the orbit of $p$ modulo $Aut(\cA)$ is included within $M^{\cE_\cA}$.
By hypothesis the number of the orbits is $k$; hence, the number of all possible values for $M^{\cE_\cA}$ cannot overcome $2^k$.

In conclusion, $A^-$ is non-empty, $\gb$-co-r.e.  and modulo $\Th(\cE_\cA)$ is finite. 
Then (iii) follows from Theorem~\ref{cor}.
\end{proof}

All the material developed in Section~\ref{section-effective}
could be adapted to the stable semantics (Berry's ccc of $DI$--domains and stable functions)
 and strongly stable semantics (Ehrhard's ccc of $DI$-domains
with coherence and strongly stable functions).
We recall that the notion of an effectively given DI-domain has been introduced by Gruchalski in \cite{Gru}, 
where it is shown that the category  having effective DI-domains as objects and stable functions as morphisms is a ccc.
There are also many effective models in the stable and strongly stable semantics. 
Indeed, the stable semantics contains a class which is analogous to the class of graph models (see \cite{Berli}),
namely Girard's class of \emph{reflexive coherent spaces} called \emph{$G$-models} in \cite{Berli}.
The results shown in Theorem~\ref{effective-graph} and in Theorem~\ref{effective-partial} 
for graph models could also be adapted 
for $G$-models, even if it is more delicate to complete partial pairs in this case 
(see \cite{Ker}).
It could also be developed for Ehrhard's class of strongly stable $H$-models (see \cite{Berli})
even though working in the strongly stable semantics certainly adds technical difficulties.

\begin{theorem}\label{Thm:8.2} Let $\cD$ be an effective $\gl$-model
in the stable or strongly stable semantics. Then $\Th(\cD)$ is not r.e. 
\end{theorem}

\begin{proof} Since $\bot_D\in \cD^{dec}$ and the interpretation function is computable, 
then   $\bot_D^- = \{ M\in \gL^o : M^\cD = \bot_D\}$ is co-r.e.
If we show that this set
is non-empty, then
$\Th(\cD)$ cannot be r.e. Since $\cD$ is effective, then by Theorem~\ref{cor}(i)
 $\ssF^-\cap \ssT^-$ is a non-empty and co-r.e. set of $\gl$-terms.
Let $N\in \ssF^-\cap \ssT^-$ and
let $f,g,h: \cD\to\cD$ be three (strongly) stable functions such that
$f(x) = \ssT^\cD \cdot x$,  $g(x) = \ssF^\cD \cdot x$ and $h(x) = N^\cD \cdot x$ for all $x\in D$. 
By monotonicity we have $h\leq_{s}f,g$ in the stable ordering. 
Now, $g$ is the constant function taking value $\ssi^\cD$, and 
$f(\bot_D)= \ssT^\cD \cdot \bot_D$. 
The first assertion forces $h$ to be a constant function, because in the stable ordering all functions under a constant map
are also constant, while the second assertion together with the fact that $h$ is pointwise smaller than $f$
 forces the constant function $h$ to satisfy $h(x) = \ssT^\cD \cdot \bot_D$ for all $x$.
Then an easy computation provides that $(NPP)^\cD = \bot_D$ for every closed term $P$.
In conclusion, we have that $\{ M\in \gL^o : M^\cD = \bot_D\}\neq \emptyset$ and the theory of $\cD$ is not r.e.
\end{proof}

\vspace*{-2.0mm}
\section{The L\"owenheim-Skolem theorem}\label{skolem}
\vspace*{-1.0mm}
In this section we show that for each graph model $\cG$ there is a countable graph model $\cP$
with the same equational/order theory. This result is a kind of downwards L\"owenheim-Skolem theorem
for graph models which positively 
 answers Question 3 in \cite[Section 6.3]{Berli}.
Note that we cannot apply directly the classical L\"owenheim-Skolem theorem since graph models
are not first-order structures.

Let $\cA,\cB$ be partial pairs.  We say that $\cA$ is a {\em subpair} of $\cB$, and we write $\cA\leq \cB$,
 if $A\subseteq B$ and  $c_\cB(a,\ga) = c_\cA(a,\ga)$ for all $(a,\ga)\in dom(c_\cA)$.

As a matter of notation, if $\gr, \gs$ are environments and $C$ is a set, we let
$\gs = \gr\cap C$ mean
$\gs(x) = \gr(x)\cap C$ for every variable $x$,
and $\gr\subseteq \gs$ mean $\gr(x)\subseteq \gs(x)$ for every variable $x$.

The proof of the following lemma is straightforward.
Recall that the definition of interpretation with respect to a partial pair is defined in Section~\ref{graph-models}.

\begin{lemma}\label{claim1.4} Suppose $\cA\leq\cB$, then
$M^\cA_\gr\subseteq M^\cB_\gs$ for all environments $\gr: Var\to \cP(A)$ and 
 $\gs: Var\to \cP(B)$ such that $\gr\subseteq\gs$.
\end{lemma}

\begin{lemma}\label{claim??} Let $M$ be a $\gl$-term, $\cG$ be a graph model and
$\ga\in M^\cG_\gr$ for some environment $\gr$.
Then there exists a finite subpair $\cA$ of $\cG$  such that
$\ga\in M^\cA_{\gr\cap A}$.
\end{lemma}

\begin{proof} The proof is by induction on $M$.

If $M\equiv x$, then $\ga\in \gr(x)$, so that we define $A = \{ \ga\}$ and $dom(c_\cA)=\emptyset$.

If $M\equiv \gl x.P$, then $\ga \equiv c_\cG(b,\gb)$ for some $b$ and $\gb$ such that $\gb\in P^\cG_{\gr[x:=b]}$.
By the induction hypothesis there exists a finite subpair $\cB$ of $\cG$ such that $\gb\in P^\cB_{\gr[x:=b]\cap B}$.
We define another finite subpair $\cA$ of $\cG$ 
as follows: $A = B\cup b\cup \{ \gb, \ga\}$;
 $dom(c_\cA) = dom(c_\cB)\cup \{ (b,\gb)\}$.
Then we have that $\cB\leq \cA$ and
$\gr[x:=b]\cap B\subseteq \gr[x:=b]\cap A$.
From $\gb\in P^\cB_{\gr[x:=b]\cap B}$ and from Lemma~\ref{claim1.4} it follows that 
 $\gb\in P^{\cA}_{\gr[x:=b]\cap A}= P^{\cA}_{(\gr\cap A)[x:=b]}$.
Then we have that $\ga \equiv c_\cA(b,\gb) \in (\gl x.P)^{\cA}_{\gr\cap A}$.

If $M\equiv PQ$, then there is $a = \{ \ga_1,\dots,\ga_n\}$ such that 
$c_\cG(a,\ga)\in P^\cG_\gr$ and $a\subseteq Q^\cG_\gr$.
By the induction hypothesis there exist finite subpairs $\cA_0, \cA_1,\dots,\cA_n$ of $\cG$ 
such that $c_\cG(a,\ga)\in P^{\cA_0}_{\gr \cap A_0}$
and $\ga_k\in Q^{\cA_k}_{\gr \cap A_k}$ for $k=1,\dots,n$.  
We define another finite subpair $\cA$ of $\cG$ as follows:
 $A = \cup_{0\leq k\leq n} A_k \cup a \cup \{\ga\}$ and
$dom(c_\cA) = (\cup_{0\leq k\leq n} dom(c_{\cA_k})) \cup \{ (a,\ga)\}$.
The conclusion follows from Lemma~\ref{claim1.4}.
\end{proof}

\begin{proposition}\label{claimclaim} 
Let $\cG$ be a graph model, and suppose $\ga\in M^\cG - N^\cG$ for some $M,N\in\gL^o$.
Then there exists a finite  $\cA\leq\cG$ such that: for all pairs $\cC\geq \cA$, if there is a morphism
$f:\cC\to\cG$ such that $f(\ga) = \ga$, then $\ga\in M^\cC - N^\cC$.
\end{proposition}
\vspace*{-1.0mm}
\begin{proof} 
By Lemma~\ref{claim??} there is a finite pair $\cA$ such that $\ga\in M^\cA$. 
By Lemma~\ref{claim1.4} we have $\ga\in M^\cC$. 
Now, if $\ga\in N^\cC$ then, by Lemma~\ref{mor} $\ga = f(\ga)\in N^\cG$, which is a contradiction.
\end{proof}

\begin{corollary}\label{claim1.7} Let $\cG$ be a graph model, and suppose $\ga\in M^\cG - N^\cG$ for some $M,N\in\gL^o$.
Then there exists a finite  $\cA\leq\cG$ such that: for all pairs $\cB$ satisfying $\cA\leq\cB\leq\cG$
we have  $\ga\in M^\cB - N^\cB$.
\end{corollary}

Let $\cG$ be a graph model. A graph model $\cP$ is called a {\em sub graph model}
of $\cG$ if $\cP\leq\cG$.
It is easy to check that the class of sub graph models of $\cG$ is closed under (finite and infinite) intersection.
If $\cA\leq\cG$ is a partial pair, 
then the {\em sub graph model generated by $\cA$}
is defined as the intersection of all graph models $\cP$ such that $\cA\leq\cP\leq\cG$.

\begin{theorem} {\rm (L\"owenheim-Skolem Theorem for graph models)} 
For every graph model $\cG$ there exists
a sub graph model $\cP$ of $\cG$ with a countable carrier set and such that $\Or(\cP) = \Or(\cG)$, and hence
$\Th(\cP) = \Th(\cG)$.
\end{theorem}

\begin{proof} 
We will define an increasing sequence of countable subpairs $\cA_n$ of $\cG$, and
take for $\cP$ the sub graph model of $\cG$ generated by
$\cA \equiv \cup \cA_n$.

First we define $\cA_0$. Let $I$ be the countable set of inequations between closed $\gl$-terms which fail in $\cG$.
Let $e\in I$. By Corollary~\ref{claim1.7} there exists a finite partial pair
$\cA_e\leq \cG$ such that $e$ fails in every partial pair $\cB$ satisfying  $\cA_e\leq \cB\leq\cG$. Then we define 
$\cA_0 = \cup_{e\in I} \cA_e \leq \cG$.
Assume now that $\cA_{n}$ has been defined.
We define $\cA_{n+1}$ as follows.
For each inequation $e \equiv M\ \sq\ N$ 
which holds in $\cG$ and fails in the sub graph model $\cP_{n}\leq \cG$ generated by $\cA_n$,
we consider the set
$L_{e} = \{ \ga \in P_{n} : \ga\in M^{\cP_n} - N^{\cP_n}\}$. Let $\ga\in L_{e}$.
Since $\cP_n\leq \cG$ and $\ga\in M^{\cP_n}$, then by Lemma~\ref{claim1.4} we have that
$\ga \in  M^{\cG}$. By $\cG\models M\ \sq\ N$ we also obtain $\ga \in N^{\cG}$.
By Lemma~\ref{claim??} there exists a partial pair $\cF_{\ga,e}\leq \cG$ 
such that $\ga\in N^{\cF_{\ga,e}}$.
We define $\cA_{n+1}$ as the union of the partial pair $\cA_n$
and the partial pairs $\cF_{\ga, e}$  for every $\ga\in L_{e}$. 

Finally take for $\cP$ the sub graph model of $\cG$ generated by $\cA \equiv \cup\cA_n$. 
By construction we have, for every inequation $e$ which fails in $\cG$:
$\cA_{e}\leq\cP_{n}\leq\cP\leq\cG$. 
Now, $\Or(\cP)\subseteq \Or(\cG)$
follows from Corollary~\ref{claim1.7} and from the choice of $\cA_e$.

Let now $M\ \sq\ N$ be an inequation which fails in $\cP$ but not in $\cG$.
Then there is  an $\ga\in M^{\cP}- N^{\cP}$. By Corollary~\ref{claim1.7}
there is a finite partial pair $\cB\leq\cP$ satisfying the following condition:
 for every partial pair $\cC$ such that $\cB\leq\cC\leq \cP$,
we have $\ga\in M^{\cC} - N^{\cC}$. 
Since $B$ is finite, we have that $\cB\leq \cP_n$ for some $n$.
This implies that $\ga\in M^{\cP_n}- N^{\cP_n}$.
By construction of $\cP_{n+1}$ we have that $\ga \in N^{\cP_{n+1}}$; this implies 
$\ga \in N^{\cP}$. Contradiction.
\end{proof}

\vspace*{-2.0mm}
\section{The minimum order graph theory}
\vspace*{-1.0mm}
In this section we show one of the main theorems of the paper: the minimum order graph theory exists and it
is the theory of an effective graph model. 
This result has the interesting consequence that no order graph theory can be r.e.

\begin{lemma}\label{claim1.3} 
Suppose $\cA\leq \cG$ and let $f:E_\cA\to G$ be defined by induction
 over the rank of $x\in E_\cA$ as follows:
$$
    f(x) = \left\{
    \begin{array}{ll}
x & \textrm{if $x\in A$}\\
c_\cG(fa,f\ga) & \textrm{if $x\notin A$ and $x\equiv (a,\ga)$.}\\
\end{array} 
    \right.
$$
Then $f$ is a morphism from $\cE_\cA$ into $\cG$.
\end{lemma}

\begin{lemma}\label{claim1.7bis} 
Suppose $\ga \in M^{\cG}-N^{\cG}$
for some $M,N\in \gL^o$. Then there exists a finite $\cA\leq \cG$ such that: for
all pairs $\cB$ satisfying 
$\cA\leq \cB\leq \cG$, we have $\ga \in M^{\cE_\cB}-N^{\cE_\cB}$.
\end{lemma}
\begin{proof} 
By Proposition~\ref{claimclaim} and Lemma~\ref{claim1.3}.
\end{proof}

\begin{theorem}\label{minimum} There exists an effective graph model
whose order/equational theory is the minimum  order/equational graph theory.
\end{theorem}

\begin{proof} 
It is not difficult to define an effective bijective numeration $\cN$ of all finite
partial pairs whose carrier set is a subset of $\nat$. We denote by $\cN_k$  the $k$-th finite partial pair
with $N_k\subseteq \nat$. We now make the carrier sets $N_k$ ($k\in\nat$) disjoint.
Let $p_k$ be the $k$-th prime natural number.
Then we define another finite partial pair $\cP_k$ as follows:
$P_k = \{p_k^{x+1}: x\in N_k\}$ and 
$c_{\cP_k}(\{p_k^{\ga_1+1},\dots,p_k^{\ga_n+1}\},p_k^{\ga +1}) = p_k^{c_{\cN_k}(\{\ga_1,\dots,\ga_n\},\ga)+1}$
for all $(\{\ga_1,\dots,\ga_n\},\ga)\in dom(c_{\cN_k})$.
In this way we get an effective bijective numeration of all finite partial pairs $\cP_k$.
Finally, we take $\cP\equiv \cup_{k\in\nat} \cP_k$.
It is an easy matter to prove that $P$ is a decidable subset of $\nat$ and that, after encoding, 
$c_\cP = \cup_{k\in\nat} c_{\cP_k}$ is a computable map
with a decidable domain and range. 
Then by Theorem~\ref{effective-partial}(ii)   $\cE_\cP$ is an effective graph model.
Notice that $\cE_\cP$ is also isomorphic to the completion of the union $\cup_{k\in\nat} \cE_{\cP_k}$,
where $\cE_{\cP_k}$ is the completion of the partial pair $\cP_k$.

We now prove that the order theory of $\cE_\cP$ is the minimum one.
Let $e\equiv M\ \sq\ N$ be an inequation which fails in some graph model $\cG$. 
By Lemma~\ref{claim1.7bis} $e$ fails in the  completion of a finite partial pair $\cA$.
Without loss of generality, we may assume that the carrier set of $\cA$ is a subset of $\nat$,
and then that $\cA$ is one of the partial pairs $\cP_k$. For such a $\cP_k$, $e$ fails in $\cE_{\cP_k}$.
Now, it was shown by Bucciarelli and Salibra in \cite[Proposition~2]{BSa} that, if a graph model $\cG$ is the completion
of the disjoint union of a family of graph models $\cG_i$, then $Q^{\cG_i} = Q^\cG \cap G_i$ for any closed $\gl$-term $Q$.
Then we can conclude the proof as follows:
if the inequation $e$ holds in $\cE_\cP$, then by \cite[Proposition~2]{BSa} we get a contradiction:
 $M^{\cE_{\cP_k}} = M^{\cE_\cP} \cap E_{P_k}\subseteq N^{\cE_\cP} \cap E_{P_k}= N^{\cE_{\cP_k}}$.
\end{proof}

\begin{theorem}\label{minimum-bis} Let $\cT_{min}$ and $\cO_{min}$ be, respectively, 
the minimum equational graph theory
and the minimum order graph theory. We have:
\begin{itemize}
\item[(i)] $\cO_{min}$ is not r.e.
\item[(ii)] $\cT_{min}$ is an intersection of a countable set of non-r.e. equational graph theories.
\end{itemize}
\end{theorem}
\begin{proof}
(i) follows from Theorem~\ref{minimum} and from Theorem~\ref{cor}(iii), because
$\cO_{min}$ is the theory of an  effective $\gl$-model.

(ii) By the proof of  Theorem~\ref{minimum} we have that $\cT_{min}$
is an intersection of a countable set of graph theories, which are theories of completions of finite partial pairs.
By Theorem~\ref{effective-partial}(iii) these theories are not r.e.
\end{proof}

\begin{corollary}\label{every} 
For all graph models $\cG$, $\Or(\cG)$ is not r.e.
\end{corollary}
\begin{proof} If $\Or(\cG)$ is r.e. and $M$ is closed and $\gb$-normal, 
 then  $M^- = \{ N\in\gL^o : N^\cG\subseteq M^\cG\}$
is a $\gb$-r.e. set,
which contains the $\gb$-co-r.e. set $\{ N \in \gL^o: \cO_{min} \vdash N\ \sq\ M\}$.
By the FIP $M^- = \gL^o$. By the arbitrariness of $M$, it follows that
$\ssT^- = \ssF^-$. Since $\ssF \in \ssT^-$ and
conversely we get $\ssF=\ssT$ in $\cG$, contradiction.
\end{proof}

\begin{corollary}
Let $\mathfrak{G}$ be the class of all graph models.
For any finite sequence $M_1,\dots, M_n$ of closed $\gb$-normal forms, 
there exists a non-empty $\gb$-closed co-r.e. set $\cU$ of closed unsolvable terms 
such that 
$$(\forall\cG\in \mathfrak{G})(\forall U\in\cU)\ U^\cG\subseteq M_1^\cG\cap\dots\cap 
M_n^\cG.$$
\end{corollary}
\begin{proof} By Theorem~\ref{cor}(i) applied to any effective graph model with minimum theory,
we have $(\forall U\in\cU)\ \cO_{min}\vdash U\ \sq\ M_1\wedge\dots\wedge \cO_{min}\vdash U\ \sq\ M_n$.
The conclusion follows.
\end{proof}
The authors do not know any example of unsolvable satisfying the above condition.

\end{document}